\documentclass[11pt]{amsart}

\usepackage[a4paper,hmargin=3.5cm,vmargin=3.5cm]{geometry}
\usepackage{amsfonts,amssymb,amscd,amstext}
\usepackage{graphicx}
\usepackage[dvips]{epsfig}

\usepackage{fancyhdr}
\usepackage{times}

\usepackage{enumerate}
\usepackage{titlesec}
\usepackage{mathrsfs}

\def\Ccal{\mathcal{C}}
\def\Bcal{\mathcal{B}}

\def\Kscr{\mathscr{K}}
\def\Lscr{\mathscr{L}}

\def\b{\mathbb{B}}

\def\r{\mathbb{R}}
\def\n{\mathbb{N}}
\def\s{\mathbb{S}}
\def\mgot{\mathfrak{m}}

\titleformat{\section}
{\filcenter\bfseries\large} {\thesection{.}}{0.2cm}{}
\titleformat{\subsection}[runin]
{\bfseries} {\thesubsection{.}}{0.15cm}{}[.]
\titleformat{\subsubsection}[runin]
{\em}{\thesubsubsection{.}}{0.15cm}{}[.]

\usepackage[up,bf]{caption}

\newtheorem{theorem}{Theorem}[section]

\newtheorem{claim}[theorem]{Claim}

\newtheorem{corollary}[theorem]{Corollary}

\newtheorem{definition}[theorem]{Definition}

\theoremstyle{definition}

\numberwithin{equation}{section}
\numberwithin{figure}{section}

\usepackage{color}

\begin{document}

\fancyhead[CO]{Capillary surfaces inside polyhedral regions} 
\fancyhead[CE]{A. Alarc\'{o}n$\,$  and$\,$ R. Souam} 
\fancyhead[RO,LE]{\thepage} 

\thispagestyle{empty}

\vspace*{1cm}
\begin{center}
{\bf\LARGE Capillary surfaces inside polyhedral regions}

\vspace*{0.5cm}

{\large\bf Antonio Alarc\'{o}n$\;$ and$\;$ Rabah Souam}
\end{center}

\footnote[0]{\vspace*{-0.4cm}

\noindent A. Alarc\'{o}n

\noindent Departamento de Geometr\'{\i}a y Topolog\'{\i}a, Universidad de Granada, E-18071 Granada, Spain.

\noindent e-mail: {\tt alarcon@ugr.es}

\vspace*{0.1cm}

\noindent R. Souam

\noindent Institut de Math\'{e}matiques de Jussieu-Paris Rive Gauche,   UMR 7586, B\^{a}timent Sophie Germain,  Case 7012, 75205  Paris Cedex 13, France.

\noindent e-mail: {\tt souam@math.jussieu.fr}

\vspace*{0.1cm}

\noindent A. Alarc\'{o}n is supported by Vicerrectorado de Pol\'{i}tica Cient\'{i}fica e Investigaci\'{o}n de la Universidad de Granada.

\noindent A.\ Alarc\'{o}n's research is partially supported by MCYT-FEDER grant MTM2011-22547 and Junta de Andaluc\'{i}a Grant P09-FQM-5088.

\noindent R.\ Souam  was partially supported by the French-Spanish Research Network in Geometric Analysis.}

\vspace*{1cm}

\begin{quote}
{\small
\noindent {\bf Abstract}\hspace*{0.1cm} In this paper we provide a large new family of embedded capillary surfaces inside polyhedral regions in the Euclidean space.   The angle of contact of the examples we furnish lies in  $[\frac{\pi}{2}, \pi]$ and it is allowed to vary from one boundary component to the other.
\vspace*{0.2cm}

\noindent{\bf Keywords}\hspace*{0.1cm} Capillary surfaces, constant mean curvature surfaces.

\vspace*{0.2cm}

\noindent{\bf Mathematics Subject Classification (2010)}\hspace*{0.1cm} 53A10, 49Q10, 53C42, 76B45.}

\end{quote}


\section{Introduction}\label{sec:intro}

Consider a (closed) region $\mathcal B$ in the Euclidean space $\r^3$. A capillary surface in $\mathcal B$ is a  compact {\em $H$-surface} (i.e.,  with constant mean curvature $H$) with non-empty boundary, which is $\Ccal^1$ up to the boundary and  meeting  the frontier $\partial \mathcal B$ of $\mathcal B$ at a constant angle $\theta \in [0,\pi]$ along its boundary.  Capillary surfaces are stationary surfaces for an energy functional under a volume constraint. More precisely, given a compact surface $\Sigma$ inside $\mathcal B$ such that $\partial \Sigma \subset \partial \mathcal B$ and $\partial \Sigma$ bounds a compact domain $W$ in
$\partial \mathcal B,$ the energy of $\Sigma$ is  by definition the quantity
\[
\mathcal E(\Sigma):= {\rm Area} (\Sigma) -\cos\theta \,  {\rm Area (W)}.
\]

The stationary surfaces of $\mathcal E$ for variations preserving the enclosed volume are precisely the 
{$H$-surfaces}  which make a constant angle $\theta$ with  $\partial \mathcal B.$ Here the contact angle is computed from inside the domain enclosed by $\Sigma\cup W.$ Capillary surfaces model incompressible liquids inside  a container in the absence of gravity; $\Sigma$ represents the free surface of the 
liquid and $W$ the region of the container wetted by the liquid. A standard reference  on capillary surfaces is the book by Finn \cite{Fi}. 

When the angle of contact $\theta$ is $\pi/2$, that is, when the surface  is orthogonal to $\partial \mathcal B$, the surface is said to be an $H$-surface with free boundary.  Many works have been devoted to prove the existence of disk-type $H$-surfaces with free boundary inside a given compact region $\Omega$  in $\r^3$. For instance, Struwe \cite{struwe1} proved the existence of a minimal disk (i.e., with $H=0$) with free boundary, which is not necessarily embedded,  in any  $\Omega$ with sufficiently regular boundary. Improving Struwe's result, Gr$\ddot{\text{u}}$ter and Jost \cite{gruter-jost} proved the existence of an embedded minimal disk with free boundary in any {\it convex}  domain $\Omega.$ For  $H\ne 0$, existence results of disk-type $H$-surfaces with free boundary for various values of $H$ have been obtained by Struwe \cite{struwe2} and Fall \cite{fall}. However no existence result of this nature is available in the literature when the angle of contact $\theta$ is not $\pi/2.$ 

In this paper we prove the existence of a large family of embedded capillary surfaces of genus zero in polyhedral convex domains in $\r^3$.  We moreover allow the angle of contact   to vary from one boundary component to the other. The angles of contact can take any value in   $[\frac{\pi}{2}, \pi]$ provided they satisfy a mild condition (see the statement of Theorem \ref{th:intro}).
 In the physical interpretation,
one allows the bounding faces  of the polyhedral container to be composed of different (homogeneous) materials. {If $F_j, j=1,\ldots,\mgot$, denote the faces of the polyhedral region $\mathcal B$, then the associated energy functional in that case is defined as follows. Given a compact surface $\Sigma$ inside $\mathcal B$ having $\mgot$ boundary components, 
$\gamma_j
\subset F_j$, bounding domains $W_j\subset F_j, j=1,\ldots,\mgot$, the energy of $\Sigma$ is
\begin{equation}\label{eq:energy}
\mathcal E(\Sigma):=  {\rm Area} (\Sigma) -\sum_{j=1}^\mgot \cos \theta_j\,  {\rm Area}(W_j),
\end{equation}
where $\theta_j$, $j=1,\ldots\mgot$, are given constants in $[0,\pi]$.
The stationary  surfaces of $\mathcal E$ for variations  preserving the enclosed volume are, this time, 
the $H$-surfaces which meet $F_j$ at a constant angle $\theta_j$ for each $j=1,\ldots, \mgot.$
 When $\Sigma$ is embedded and either $\mgot=1$, or $\mgot =2$ and $\mathcal B$ is a slab,  then, using the Alexandrov reflection technique, one shows $\Sigma$ is a rotational surface and so it is a part of a sphere, a part of a cylinder, or a part of  a Delaunay surface (either an unduloid or a nodoid). Wente \cite{We}, has constructed an immersed, non embedded, annulus of constant mean curvature intersecting orthogonally two parallel planes.
 In case $\mathcal B$ is a wedge, that is, the region of the space between two intersecting planes, Park \cite{Pa} has shown that if $\Sigma$ is  embedded, topologically an annulus,  and does not touch the vertex of the wedge, then it is part of a round sphere (see also  McCuan \cite{Mc} for a partial result).  }

Recently the authors \cite{AS1} found polyhedral regions in $\r^3$ admitting capillary surfaces with contact angle equal to $\pi$ along every boundary component. In this paper we provide a large new family of capillary surfaces inside polyhedral regions, widely generalizing our previous results. 
The following is our main result.

\begin{theorem}\label{th:intro}
Let $\{p_1,\ldots,p_\mgot\}$ be a subset of  $\s^2$ with cardinal number $\mgot\in\n$ and let $H>0$ be a positive real number. Let also  $(\theta_1,\ldots,\theta_\mgot)\in(\frac{\pi}{2}, \pi)^\mgot$ be given numbers. 
Denote by  $B_j$ the open disk in $\s^2$ of radius $\pi-\theta_j$ centered at $p_j, j=1,\ldots,\mgot$, and suppose the closed disks $\overline B_j$ are pairwise disjoint. 

Consider real numbers $a_j>0$, $j=1,\ldots,\mgot$, such that the balancing condition
\begin{equation}\label{eq:balancing}
\sum_{j=1}^\mgot \Big(a_j-\frac{\pi}{4H^2} \sin^2\theta_j \Big) p_j =\vec{0}\in\r^3
\end{equation}
is satisfied.

Then there exists a polyhedral region $\mathcal B$ in $\r^3$ with frontier $\partial \mathcal B$ consisting of $\mgot$ planar regions $F_1,\ldots,F_\mgot$, where $F_j$ is orthogonal to $p_j$ for all $j\in\{1,\ldots,\mgot\}$, and there exists an embedded capillary surface $\Sigma$ in $\mathcal B$, with constant mean curvature $H$, satisfying the following properties:
\begin{enumerate}[\rm (i)]
\item $\Sigma$ is $\Ccal^{1,1}$ up to the boundary.
\item $\gamma_j:=\Sigma\cap F_j$ is a convex Jordan curve contained in the relative interior of $F_j$. 
\item $\Sigma$ meets $F_j$ at constant angle $\theta_j$ along $\gamma_j$ for all $j\in\{1,\ldots,\mgot\}$.
\item $\Sigma\setminus\partial\Sigma$ is positively curved and $\Sigma\cup (\cup_{j=1}^\mgot D_j)$ is the boundary surface of a convex body in $\r^3,$ where $D_j$ is the convex connected component of  $F_j\setminus\gamma_j$ for all $j\in\{1,\ldots,\mgot\}$.
\item $a_j= {\rm Area}(D_j)-{\rm Length}(\partial D_j)\frac{\sin\theta_j}{2H} + \pi \left(\frac{\sin\theta_j}{2H}\right)^2$  for all $ j\in\{1,\ldots,\mgot\}.$
\end{enumerate}
Furthermore, the data $(\mathcal B, \Sigma)$ satisfying the previous conditions are unique up to translations.

Conversely 
let $\{p_1,\ldots,p_\mgot\}\subset\s^2$, let $H>0$, let $(\theta_1,\ldots,\theta_\mgot)\in(\frac{\pi}{2}, \pi)^\mgot$, and let $B_j\subset\s^2$, $j=1,\ldots,\mgot$, be as above. Assume there exists a polyhedral region $\Bcal$ in $\r^3$ with frontier $\partial\Bcal$ consisting of $\mgot$ planar regions $F_1,\ldots, F_\mgot$, where $F_j$ is orthogonal to $p_j$ for all $j=1,\ldots,\mgot$, and there exists a capillary surface $\Sigma$ which is $\mathcal C^1$ up to the boundary, with constant mean curvature $H$, in $\Bcal$ enjoying properties  {\rm (ii)}--{\rm (iv)}. 
Suppose in addition that $\Sigma$ is not a part of an unduloid.

Then the numbers $a_j$, $j=1,\ldots,\mgot$, defined by {\rm (v)} are all positive and satisfy the balancing condition \eqref{eq:balancing}.
\end{theorem}

Consider a vertical unduloid $S\subset\r^3$, of mean curvature $H>0$, and a horizontal plane $\Pi\subset\r^3$ such that $S$ meets $\Pi$ at a constant angle $\theta\in(\frac{\pi}2,\pi)$. Obviously, the intersection $S\cap\Pi$ bounds a round circle, $D\subset\Pi$, and one can check that for a suitable value of $\theta$ it holds that ${\rm Area}(D)-{\rm Length}(\partial D)\frac{\sin\theta}{2H} + \pi \left(\frac{\sin\theta}{2H}\right)^2=0$.  In particular, Theorem \ref{th:intro} does not apply to pieces of unduloids having such a particular boundary component.

Our main tool in the proof of the existence part of Theorem \ref{th:intro} is {\em the classical Minkowski problem of prescribing positive Gauss curvature on the sphere $\s^2$}; this will be done in Sec.\ \ref{sec:proof}. We solve such problem for a sequence of curvature functions on $\s^2$ in order to construct a suitable sequence of strictly convex bodies in $\r^3$. We then find the surface $\Sigma$ satisfying the conclusion of Theorem \ref{th:intro} as a parallel surface to a piece of the boundary of the convex body 
limit of that sequence. On the other hand, the uniqueness part of Theorem \ref{th:intro} and the fact that the balancing condition \eqref{eq:balancing} is in fact necessary are proved in Sec.\ \ref{sec:balancing}.

Replacing condition {\rm (i)} by the weaker one $\Sigma$ is $\Ccal^1$ up to the boundary, then the statement of Theorem \ref{th:intro} in the limit case $\theta_j=\pi$ and $B_j=\{p_j\}$ for all $j=1,\ldots,\mgot$ (hence $\mgot\geq 2$ by \eqref{eq:balancing}) was recently proven by the authors in \cite{AS1}. Combining the arguments in the present paper with those in \cite{AS1}, it can be also showed that Theorem \ref{th:intro} remains valid if one allows some of the angles $\theta_j$ to be equal to $\pi$; however only $\Ccal^1$ regularity is ensured in the boundary components of $\Sigma$  corresponding to those angles. It remains open the question whether these examples and those in \cite{AS1} are $\Ccal^{1,1}$ up to the boundary. Also it would be interesting to decide whether the capillary surfaces furnished by Theorem \ref{th:intro} enjoy higher boundary regularity; we briefly discuss this issue in Sec.\ \ref{sec:comments}.

An important remark is that symmetries in the initial data of Theorem \ref{th:intro} are induced to the resulting capillary surfaces. More specifically, if $s\colon\r^3\to\r^3$ denotes the symmetry with respect to a vectorial plane $\Pi\subset\r^3$ and there exists a bijective map $\eta\colon\{1,\ldots,\mgot\}\to\{1,\ldots\mgot\}$ such that $B_{\eta(j)}=s(B_j)$ and $a_{\eta(j)}=a_j$ for all $j$, then, up to a translation, the polyhedral region $\Bcal$ and the capillary surface $\Sigma$ furnished by Theorem \ref{th:intro} are symmetric with respect to $\Pi$. Indeed, just observe that the uniqueness part of Theorem \ref{th:intro} ensures that, up to a translation, $s(\Bcal)=\Bcal$ and $s(\Sigma)=\Sigma$. As application we obtain in Sec.\ \ref{sec:pi2} the existence of capillary surfaces as those in Theorem \ref{th:intro} in which one of the contact angles is allowed to be $\frac{\pi}2$ (since the closed disks $\overline B_1,\ldots, \overline B_\mgot$ in Theorem \ref{th:intro} are assumed to be pairwise disjoint, then at most one of the angles $\theta_j$ could be equal to $\frac{\pi}2$).

\begin{corollary}\label{co:pi2}
The statement of Theorem \ref{th:intro} remains true if one of the contact angles equals $\frac{\pi}2$.
\end{corollary}

Let $\s^2(r)\subset\r^3$ denote the sphere of radius $r>0$ centered at $\vec 0$. Cutting $\s^2(1/H)$ by a finite family of affine planes one gets the surfaces corresponding to the cases where $a_j= \frac{\pi}{4H^2}\sin^2\theta_j$ for all $j=1,\ldots,\mgot.$ For instance, if $\mgot=1$ then \eqref{eq:balancing} implies $a_1=\frac{\pi}{4H^2}\sin^2\theta_1$ and $\Sigma\subset \s^2(1/H)$. Likewise, if $\mgot=2$ and $p_1$ and $p_2$ are not antipodal, then necessarily $a_j=\frac{\pi}{4H^2} \sin^2\theta_j, \, j=1,2$, and the surfaces we obtain are again pieces of the sphere $\s^2(1/H).$
On the other hand, if $\mgot = 2$ and $p_2=-p_1$, then the solutions are surfaces of revolution as is shown by the Alexandrov reflection technique; recall that $\Sigma$ is $\Ccal^1$ up to the boundary (see {\rm (i)}).

Suppose now $\mgot \geq 3$  and let  $H>0$ be a real positive number.  Consider  $(\theta_1,\ldots,\theta_{\mgot}) \in (\frac{\pi}{2}, \pi)^{\mgot}$, and 
$\{p_1,\ldots, p_\mgot\} \subset \s^2$ satisfying the hypotheses of Theorem \ref{th:intro}. Denote by $\mathcal F_H(\theta_1,\ldots,\theta_\mgot;  p_1,\ldots,p_\mgot)$ the family of capillary surfaces, up to translations, provided by Theorem \ref{th:intro}  for all possible values of  $a_1,\ldots, a_{\mgot}$. 

If  $p_1,\ldots,p_\mgot $ are contained in an equator of $\s^2$ then $\mathcal F_H(\theta_1,\ldots,\theta_\mgot;  p_1,\ldots,p_\mgot)$ is an $(\mgot-2)$-parameter family. We may assume  without loss of generality that the points $p_1,\ldots,p_\mgot \in \s^2$ are contained in the horizontal equator, the polyhedra provided by Theorem \ref{th:intro} are then non-compact and have vertical faces. The Alexandrov reflection technique shows in this case that all the  surfaces 
in $\mathcal F_H (\theta_1,\ldots,\theta_\mgot;  p_1,\ldots,p_\mgot)$ have a horizontal plane of symmetry. 

If now $p_1,\ldots,p_\mgot $ are not contained in an equator then $\mathcal F_H(\theta_1,\ldots,\theta_\mgot;  p_1,\ldots,p_\mgot)$  is an $(\mgot-3)$-parameter family. In particular, for  $\mgot=3$, this family  reduces to a single element, namely a suitable piece of $\s^2(1/H)$.

We also point out that the Gauss map of the open $H$-surface $\Sigma\setminus\partial\Sigma$ in Theorem \ref{th:intro} is a diffeomorphism onto $\s^2\setminus\cup_{j=1}^{\mgot} \overline B_j$; see {\rm (iii)} and {\rm (iv)}. Since the Gauss map of $H$-surfaces is harmonic for the conformal structure induced by isothermal charts \cite{Ruh}, this provides harmonic diffeomorphisms $\Sigma\setminus\partial\Sigma \to \s^2\setminus\cup_{j=1}^{\mgot} \overline B_j$. The problem of determining whether there exist harmonic diffeomorphisms between given Riemannian surfaces is an important question with large literature; see e.g.\ \cite{SY,Ma,CR,AS1,AS2}.



\section{Preliminaries}\label{sec:prelim}

We denote by $\langle\cdot,\cdot\rangle$ and $\|\cdot\|$ the Euclidean inner product and norm in $\r^n$, $n\in\n$. Likewise we denote by ${\rm Vol}(\cdot)$, ${\rm Area}(\cdot)$, and ${\rm Length}(\cdot)$ the volume, area, and length operators. Given a subset $S\subset\r^n$, we write $\overline S$ for the closure of $S$ in $\r^n$ and $\partial S=\overline S\setminus S$ for its topological frontier. We call $\r_+$ the set of positive real numbers.

A topological surface is said to be {\em open} if it is non-compact and does not contain any boundary point.

A {\em convex body} in $\r^n$ is a compact convex subset of $\r^n$ with non-empty interior. A {\em strictly convex body} is a convex body whose boundary does not contain any nontrivial line segment. A compact surface in $\r^3$, with empty boundary, is said to be {\em (strictly) convex} if it is the boundary surface of a (strictly) convex body in $\r^3$.

The main tool that we exploit in this paper is {\em the classical Minkowski problem of prescribing positive Gauss curvature on the sphere}, which we now describe. 

As usual, we denote by $\s^2\subset\r^3$ the round sphere of radius $1$ centered at $\vec 0$. Let $X\colon \s^2\to\r^3$ be a $\Ccal^2$ immersion whose image surface $X(\s^2)$ is a strictly convex compact surface in $\r^3$. It trivially follows that the Gauss map $N_X\colon \s^2\to \s^2$ of $X$ is a homeomorphism and the Gauss curvature function $K\colon \s^2\to\r$ of $X$ is positive everywhere on $\s^2$. Label $\kappa\colon \s^2\to\r_+$, $\kappa=K\circ N_X^{-1}$. As observed by Minkowski, the balancing condition
\begin{equation}\label{eq:minkowski}
\int_{\s^2} \frac{p}{\kappa(p)}\, dp=0
\end{equation}
is satisfied. Obviously, the reparameterization $X\circ N_X^{-1}\colon \s^2\to \r^3$ is an embedding with Gauss curvature function $\kappa$ and Gauss map the identity map of $\s^2$.

The converse of this problem is a classical question in Differential Geometry, known in the literature as the Minkowski problem. It turns out that condition \eqref{eq:minkowski} is necessary and sufficient. The following general result is the key of our construction.

\begin{theorem}[\text{\cite{N,CY}}]\label{th:minkowski}
Let $\kappa\colon \s^2\to\r_+$ be a smooth ($\Ccal^\infty$) function satisfying the balancing condition \eqref{eq:minkowski}.

Then there exists a unique up to translations smooth embedding $X\colon \s^2\to\r^3$ such that $X(\s^2)$ is a strictly convex compact surface, the Gauss curvature function of $X$ equals $\kappa$, and the Gauss map of $X$ equals the identity map of $\s^2$.
\end{theorem}

Theorem \ref{th:minkowski} becomes useful for our purposes when joined together with the fact that {\em parallel surfaces to those with positive constant Gauss curvature in $\r^3$ are $H$-surfaces} in the following sense. Let $S\subset \r^3$ be a smooth open surface with constant Gauss curvature $K=1$; from now on, a $K$-surface. Denote by $N_S\colon S\to\s^2\subset\r^3$ the {\em outer Gauss map} of $S$; the one that at every point $p\in S$ points to the connected component of $\r^3\setminus T_pS$ disjoint from an open neighborhood of $p$ in $S$. Here we denote by $T_pS$ the affine tangent plane to the surface $S$ at the point $p$ and are taking into account that, since $K>0$, $S$ is locally strictly convex. In this setting, the smooth surface
\[
S+N_S:=\{p+N_S(p)\colon p\in S\}
\]
is locally strictly convex, has constant mean curvature $H=1/2$, and its outer Gauss map $N_{S+N_S}\colon S+N_S\to\s^2$ satisfies $N_{S+N_S}(q)=N_S(p)$ for any $q=p+N_S(p)\in S+N_S$. The surface $S+N_S$ is said to be the {\em outer parallel surface} to $S$ at distance $1$.

To finish this section let us state some notation that will be useful later on. For $q\in\s^2$ and $r\in(0,1)$, we denote by $B(q;r)$ the spherical cup of $\s^2$ given by
\[
B(q;r)=\{p\in\s^2\colon \sin\sphericalangle(p,q)<r,\; \cos\sphericalangle(p,q)>0\},
\]
where $\sphericalangle(p,q)\in[0,\pi]$ denotes the spherical angle between $p$ and $q.$ Observe that $B(q;r)$ is an open geodesic disk centered at $q$. We write $\overline B(q;r)$ for the closure of $B(q;r)$ in $\s^2$. For $s\in(0,1-r)$, we denote
\[
A(q;r,s)=B(q;r+s)\setminus\overline B(q;s);
\]
an open geodesic annulus. Likewise  we write $\overline A(q;r,s)$ for the closure of $A(q;r,s)$ in $\s^2$.

Given $q\in\s^2$ and $r\in (0,1)$, straightforward computations give
\begin{equation}\label{eq:areaBj}
{\rm Area}(B(q;r))=2\pi\big( 1-\sqrt{1-r^2} \big)
\end{equation}
and
\begin{equation}\label{eq:spherical1}
\int_{B(q;r)} p \, dp =\pi r^2 q.
\end{equation}
In particular, if $s\in(0,1-r)$, then
\begin{equation}\label{eq:spherical2}
\int_{A(q;r,s)} p\, dp=\pi(r^2+2rs)q.
\end{equation}


{\section{Existence}\label{sec:proof}

In this section we prove the existence part in Theorem \ref{th:intro} in the particular case $H=1/2$. The general case directly follows by rescaling.}

Let $\mgot\in\n$, $\{p_1,\ldots,p_\mgot\}\subset\s^2$, and $(\theta_1,\ldots,\theta_\mgot)\subset (\frac{\pi}2,\pi)^\mgot$ be as in the statement of Theorem \ref{th:intro}. Denote by  $B_j$ the open disk in $\s^2$ of radius $\pi-\theta_j$ centered at $p_j$, $j=1,\ldots,\mgot$. Finally, choose numbers $(a_1,\ldots,a_\mgot)\in(\r_+)^\mgot$ satisfying the balancing condition \eqref{eq:balancing} for $H=1/2$, that is
\begin{equation}\label{eq:balancing1/2}
\sum_{j=1}^\mgot (a_j-\pi\sin^2\theta_j)p_j=\vec 0.
\end{equation}

The first step in the proof consists of using Theorem \ref{th:minkowski} (the Minskowski problem) in order to construct for any large enough $n\in\n$ a strictly convex body $\Lscr_n$ in $\r^3$.

To begin, set 
\begin{equation}\label{eq:rj}
r_j:=\sin \theta_j\in(0,1)
\end{equation}
and fix a natural number $n_0\in\n$ large enough so that the sets $\overline B(p_j;r_j+\frac1{n_0})$, $j=1,\ldots,\mgot$, are pairwise disjoint and
\begin{equation}\label{eq:n0}
r_j+\frac1{n_0}<1 \quad \forall j\in\{1,\ldots,\mgot\}.
\end{equation}
Observe that 
\begin{equation}\label{eq:Bj}
B_j=B(p_j;r_j)\quad \forall j\in\{1,\ldots,\mgot\}.
\end{equation}
See Sec.\ \ref{sec:prelim} for notation. For any $n\geq n_0$, denote 
\begin{equation}\label{eq:Omegan}
\Delta_n:=\s^2\setminus\bigcup_{j=1}^\mgot B \Big(p_j;r_j+\frac1{n}\Big).
\end{equation}
Notice that
\begin{equation}\label{eq:s2}
\s^2=\Delta_n\cup \bigcup_{j=1}^\mgot B\Big(p_j;\frac1{n}\Big) \cup \bigcup_{j=1}^\mgot \overline A\Big(p_j;r_j,\frac1{n}\Big)
\end{equation}
and the sets $\Delta_n$, $B(p_j;\frac1{n})$, $j=1,\ldots,\mgot$, and $A(p_j;r_j,\frac1{n})$, $j=1,\ldots,\mgot$, are pairwise disjoint.

Choose a smooth ($\Ccal^\infty$) function $\kappa_n\colon\s^2\to \r_+$ satisfying
\begin{eqnarray}
\kappa_n|_{\Delta_n} & = & 1,\label{eq:f1}
\\
\kappa_n|_{B(p_j;\frac1{n})} & = & \frac{\pi}{a_j n^2},\quad \forall j=1,\ldots,\mgot, \label{eq:f11}
\\
\int_{A(p_j;r_j,\frac1{n})} \frac{p}{\kappa_n(p)}\, dp & = & \pi\left(\frac1{n^2}+\frac{2r_j}{n}\right) p_j,\quad \forall j=1,\ldots,\mgot.\label{eq:f2}
\end{eqnarray}
Such a function always exists in view of \eqref{eq:spherical2}; for instance, one can take $\kappa_n$ to be radial over each annulus $A(p_j;r_j,\frac1{n})$, $j=1,\ldots,\mgot$. We point out that $\kappa_n$ may take both large and small values on $A(p_j;r_j,\frac1{n})$, $j=1,\ldots,\mgot$.

Let us check that the Minkowski problem can be solved for the function $\kappa_n$, $n\geq n_0$. Indeed, \eqref{eq:s2} ensures that 
\begin{multline}\label{eq:minko}
\int_{\s^2} \frac{p}{\kappa_n(p)}\, dp =
\\
\int_{\Delta_n} \frac{p}{\kappa_n(p)}\, dp + \sum_{j=1}^\mgot \int_{A(p_j;r_j,\frac1{n})} \frac{p}{\kappa_n(p)}\, dp +\sum_{j=1}^\mgot \int_{B(p_j;\frac1{n})} \frac{p}{\kappa_n(p)}\, dp,
\end{multline}
recall that $\Delta_n$, $B(p_j;\frac1{n})$, $j=1,\ldots,\mgot$, and $A(p_j;r_j,\frac1{n})$, $j=1,\ldots,\mgot$, are pairwise disjoint.
Let us compute the addends in the right term of the above equation. On the one hand,
\[
\int_{\Delta_n} \frac{p}{\kappa_n(p)}\, dp \stackrel{\text{\eqref{eq:f1}}}{=} 
\int_{\Delta_n} p\, dp \stackrel{\text{\eqref{eq:Omegan}}}{=} - \sum_{j=1}^\mgot \int_{B(p_j;r_j+\frac1{n})} p\, dp,
\]
where for the latter equality we have also used that $\int_{\s^2} p\, dp=\vec 0$. Taking into account \eqref{eq:n0} and \eqref{eq:spherical1}, the above equation reads
\begin{equation}\label{eq:minko1}
\int_{\Delta_n} \frac{p}{\kappa_n(p)}\, dp =-\sum_{j=1}^\mgot\pi\left(r_j+\frac1{n}\right)^2 p_j.
\end{equation}
On the other hand, for any $j\in\{1,\ldots,\mgot\}$,
\begin{equation}\label{eq:minko2}
\int_{B(p_j;\frac1{n})} \frac{p}{\kappa_n(p)}\, dp  \stackrel{\text{\eqref{eq:f11}}}{=} 
\int_{B(p_j;\frac1{n})} \frac{a_j n^2}{\pi} p\, dp \stackrel{\text{\eqref{eq:n0}, \eqref{eq:spherical1}}}{=} a_j p_j.
\end{equation}
Combining \eqref{eq:minko}, \eqref{eq:f2}, \eqref{eq:minko1}, and \eqref{eq:minko2}, one infers that
\[
\int_{\s^2} \frac{p}{\kappa_n(p)}\, dp = 
\sum_{j=1}^\mgot (a_j-\pi r_j^2) p_j = \vec 0,
\]
where the latter equality follows from \eqref{eq:rj} and \eqref{eq:balancing1/2}. Therefore, given $n\geq n_0$, Theorem \ref{th:minkowski} can be applied to the smooth function $\kappa_n$, furnishing a smooth embedding $X_n\colon\s^2\to\r^3$ such that 
\[
S_n:=X_n(\s^2)
\]
is a strictly convex compact surface, the Gauss curvature function of $X_n$ equals $\kappa_n$, and the Gauss map of $X_n$ equals the identity map of $\s^2$. In particular, 
\begin{equation}\label{eq:dA}
\text{the area element of $X_n$ equals $\frac1{\kappa_n}$ times the one of $\s^2$,}
\end{equation}
and, in view of \eqref{eq:f1}, 
\[
\Omega_n:=X_n(\Delta_n)\subset S_n
\]
is a smooth compact $K$-surface (that is, with constant Gauss curvature $K=1$) with boundary.
Set
\[
M_{n,j}:=X_n\Big(B\Big(p_j;r_j+\frac1{n}\Big)\Big),\quad O_{n,j}:=X_n\Big(B\Big(p_j;\frac1{n}\Big)\Big),
\]
$j=1,\ldots,\mgot$, and observe that
\begin{equation}\label{eq:Sn}
S_n=\Omega_n\cup \bigcup_{j=1}^\mgot M_{n,j} = \Omega_n\cup \bigcup_{j=1}^\mgot O_{n,j} \cup \bigcup_{j=1} X_n\Big(\overline A\Big(p_j;r_j,\frac1{n}\Big)\Big)
\end{equation}
and the relative interiors of the sets involved in the latter expression are pairwise disjoint
(see \eqref{eq:s2}). Finally, denote by $\Lscr_n$ the strictly convex body bordered by $S_n$; hence
\[
S_n=\partial\Lscr_n\quad\forall n\geq n_0.
\]

The next step in the proof consists of proving that
\begin{claim}\label{cla:bounds}
There exist a uniform upper bound of the (extrinsic) diameters and a uniform lower bound of the inner diameters of $\Lscr_n$, $n\geq n_0$.
\end{claim}

This will ensure that the sequence of convex bodies $\{\Lscr_n\}_{n\geq n_0}$ has a subsequence converging in the Hausdorff distance to a convex body $\Lscr$ in $\r^3$.

Indeed, denote by $\ell_n$ the extrinsic diameter of $\Lscr_n$ and, up to a translation, assume that $\vec{0}$ is the middle point between two points $x_n$ and $y_n$ in $S_n$ at distance $\ell_n$. The following argument is a slight modification of the proof of Lemma 4.2 in \cite{AS1} (see also \cite{CY}). Let us first find $\tau>0$ such that
\begin{equation}\label{eq:arriba}
\Lscr_n\subset \b(\vec 0, \tau)\quad\forall n\geq n_0,
\end{equation}
where we denote by $\b(y,r)$ the Euclidean ball in $\r^3$ of radius $r>0$ centered at $y\in\r^3$. Set $u_n:=\frac{x_n-y_n}{\|x_n-y_n\|}$ and let
\[
h_n\colon\s^2\to\r,\quad h_n(p)=\langle p\,,\, X_n(N_{X_n}^{-1}(p))\rangle,
\]
be the support function of $ S_n$. Notice that
\[
h_n(p)=\sup_{q\in S_n} \langle p,q\rangle \geq \frac{\ell_n}2 \max\{ 0\,,\, \langle p, u_n\rangle\} \quad \forall p\in\s^2;
\] 
hence
\begin{equation}\label{eq:bound1}
\int_{\s^2} \frac{h_n(p)}{\kappa_n(p)}\, dp \geq \frac{\ell_n}{2}\int_{\s^2} \frac{\max \{ 0\,,\, \langle p, u_n\rangle\}}{\kappa_n(p)}\, dp \geq c_0 \ell_n,
\end{equation}
where 
\begin{equation}\label{eq:c0}
c_0:=2\inf_{w\in\s^2}  \int_{\Delta_{n_0}} \max\{ 0\,,\, \langle p, w\rangle\}\, dp.
\end{equation}
Here we have taken into account that $\kappa_n|_{\Delta_n}=1$ (see \eqref{eq:f1}). Observe that $c_0>0$ by \eqref{eq:n0} and does not depend on $n$.
In view of  \eqref{eq:dA} and the isoperimetric inequality in $\r^3$, one has
\begin{equation}\label{eq:bound2}
\int_{\s^2} \frac{h_n(p)}{\kappa_n(p)}\, dp = \int_{ S_n} h_n(X_n^{-1}(p))\, dp=3{\rm Vol}( \Lscr_n) \leq c_1 \sqrt{{\rm Area}( S_n)^3}
\end{equation} 
for some constant $c_1>0$ not depending on $n$. Let us find a uniform upper bound of ${\rm Area}( S_n)$, $n\geq n_0$. Indeed, from \eqref{eq:dA} and \eqref{eq:f1} we obtain that
\begin{equation}\label{eq:areaSigman}
{\rm Area}( \Omega_n) ={\rm Area}(\Delta_n).
\end{equation}
On the other hand, combining \eqref{eq:dA}, \eqref{eq:f11}, and \eqref{eq:areaBj}, one gets
\begin{equation}\label{eq:area}
{\rm Area}( O_{n,j})= 2 \big(n^2-n\sqrt{n^2-1}\big)\, a_j\quad \forall j\in\{1,\ldots,\mgot\}.
\end{equation}
Finally, let $\delta_0>0$ be a large enough constant so that $\langle p,p_j\rangle>1/\delta_0$ for all $p\in A(p_j;r_j,\frac1{n})$ and all $j\in\{1,\ldots,\mgot\}$. Taking into account \eqref{eq:dA}, \eqref{eq:f2}, and that $\|p_j\|=1$, one has
\begin{multline}\label{eq:areaAnillo}
{\rm Area}\Big(X_n\Big(A\Big(p_j;r_j,\frac1{n}\Big)\Big)\Big)<
\\
\delta_0\int_{A(p_j;r_j,\frac1{n})} \frac1{\kappa_n(p)}\langle p,p_j\rangle\, dp=
\pi\Big(\frac1{n^2}+\frac{2r_j}{n}\Big)\, \delta_0 \quad \forall j\in\{1,\ldots,\mgot\}.
\end{multline}

Joining together \eqref{eq:Sn}, \eqref{eq:areaSigman}, \eqref{eq:area}, and \eqref{eq:areaAnillo}, one has
\begin{equation}\label{eq:bound4}
{\rm Area}( S_n)<4\pi+2\sum_{j=1}^\mgot a_j+\pi\delta_0\sum_{j=1}^\mgot (1+2r_j)<c_2
\end{equation}
for some constant $c_2>0$ not depending on $n$. Combining \eqref{eq:bound1}, \eqref{eq:bound2}, and \eqref{eq:bound4} we obtain that
$\ell_n<\frac{c_1 \sqrt{c_2^3}}{c_0}$ for all $n\geq n_0$, and so \eqref{eq:arriba} holds for any $\tau>\frac{c_1 \sqrt{c_2^3}}{2c_0}$.

Let us now show the existence of $z\in \r^3$ and $\rho>0$ such that
\begin{equation}\label{eq:abajo}
\b\big(z,\rho\big)\subset \Lscr_n \quad\forall n\geq n_0.
\end{equation}
Indeed, reason by contradiction and assume that for every $k\in\n$ there exists $n_k\in\n$, $n_k\geq n_0$, such that $ \Lscr_{n_k}$ contains no ball of radius $\frac1{k}$. In this setting, \eqref{eq:arriba} and Blaschke selection theorem \cite{Sc} ensure that, up to passing to a subsequence, $\{\Lscr_n\}_{n\geq n_0}$ converges in the Hausdorff distance to a convex set $\Lscr_\infty\subset \b(\vec{0},\tau)$ which contains no Euclidean ball; hence $\Lscr_\infty$ is contained in an affine plane $\Pi_\infty\subset\r^3$. On the other hand, \eqref{eq:bound1} and \eqref{eq:bound2} give that
${\rm Area}( S_n(\Pi))\geq \frac{c_0}3$ for all $n\in\n$ and all affine plane $\Pi\subset\r^3$,
where $ S_n(\Pi)$ denotes the vertical projection of $ S_n$ on $\Pi$ (we have used here that ${\rm Vol}( \Lscr_n)\leq{\rm Area}( S_n(\Pi))\ell_n$); see \eqref{eq:c0}. However, this is in contradiction with the fact that the sequence $\{{\rm Area}( S_n(\Pi))\}_{n\geq n_0}$ converges to zero for all plane $\Pi$ orthogonal to $\Pi_\infty$. This contradiction furnishes $\rho>0$ such that for every $n\geq n_0,$ there exists a point $z_n \in \r^3$  so that $\b(z_n, 2\rho)\subset \Lscr_n$. Up to passing to a subsequence, $\{z_n\}_{n\geq n_0}$ converges to a point $z\in \r^3$ (take into account \eqref{eq:arriba}) which together with $\rho$ satisfy \eqref{eq:abajo}.

Equations \eqref{eq:arriba} and \eqref{eq:abajo} imply Claim \ref{cla:bounds}. 

In view of Claim \ref{cla:bounds}  Blaschke selection theorem \cite{Sc} applies, ensuring that, up to passing to a subsequence, $\{ \Lscr_n\}_{n\geq n_0}$ converges in the Hausdorff distance to a convex body $ \Lscr$ in $\r^3$ which consists of the accumulation set of points in $\{\Lscr_n\}_{n\geq n_0}$. 

Denote by $\widetilde S_n$ and $\widetilde \Lscr$ the outer parallel surface to $S_n$ and the outer parallel convex body to $\Lscr$ at distance $1$. Also set $S:=\partial \Lscr$ and $\widetilde S:=\partial \widetilde\Lscr$. 

The final step in the proof consists of finding a piece $\Sigma$ of $\widetilde S$ satisfying the conclusion of Theorem \ref{th:intro}.

Denote by $\Omega\subset S$ and $\overline O_j\subset S$ the compact accumulation sets of the sequences $\{\Omega_n\}_{n\geq n_0}$ and $\{O_{n,j}\}_{n\geq n_0}$, $j=1,\ldots,\mgot$. 
From \eqref{eq:areaAnillo}, the sequence $\big\{{\rm Area}\big(X_n\big(A\big(p_j;r_j,\frac1{n}\big)\big)\big)\big\}_{n\geq n_0}$ converges to $0$; hence 
$S = \Omega \cup (\cup_{j=1}^\mgot \overline O_j)$.
Further, $(\Omega\setminus\partial\Omega)\cap O_j=\emptyset$  where $O_j:=\overline O_j\setminus \partial\overline O_j$, for all $j=1,\ldots,\mgot$. In particular, 
\begin{equation}\label{eq:SUCj}
\partial \Omega =\bigcup_{j=1}^\mgot\partial O_j \quad\text{and}\quad  S= \Omega\cup \bigcup_{j=1}^\mgot  O_j.
\end{equation}

Let us check that
\begin{claim}\label{cla:discos}
$O_j$ is an open disk contained in an affine plane $\Pi_j\subset\r^3$ orthogonal to $p_j$ for all $j=1,\ldots,\mgot$.
\end{claim}

Indeed, let $N_{S_n}\colon S_n\to\s^2$ be the outer Gauss map of $S_n$ for all $n\geq n_0$. Choose $j\in\{1,\ldots,\mgot\}$ and for every $n\geq n_0$ set 
\[
C_{n,j}:=O_{n,j}+n\sqrt{\frac{a_j}{\pi}}\, \left(N_{S_n}|_{O_{n,j}}-p_j\right)
\]
the outer parallel surface to $O_{n,j}$ at distance $n\sqrt{\frac{a_j}{\pi}}$ translated by the vector $-n\sqrt{\frac{a_j}{\pi}}p_j$. 
By the very definition of $O_{n,j}$ and the fact that the Gauss map of $X_n$ is the identity map, it follows that 
\[
\|N_{S_n}(p)-p_j\|<\frac2{n}\quad\forall p\in O_{n,j},\text{ for any large enough $n\in\n$.}
\]
Therefore, we may assume that $\big\{n\sqrt{\frac{a_j}{\pi}}\big(N_{S_n}|_{O_{n,j}}-p_j\big)\big\}_{n\geq n_0}$ converges in the Hausdorff distance to a vector $v_j\in\r^3$. It follows that $v_j+\overline O_j$ is the accumulation set of the sequence $\{C_{n,j}\}_{n\geq n_0}$. 

On the other hand, \eqref{eq:f11} ensures that $O_{n,j}$ is an open surface with constant Gauss curvature equal to $\frac{\pi}{a_j n^2}$. Therefore $C_{n,j}$ is an open surface with constant mean curvature equal to $\frac1{2n}\sqrt{\frac{\pi}{a_j}}$ and positive Gauss curvature. In particular, the square of the norm of the second fundamental form of $C_{n,j}$ is bounded from above by $\frac{\pi}{4a_j n^2}$ and standard compactness arguments in Constant Mean Curvature Surface Theory (see for instance Proposition 2.3 and Lemma 2.4 in \cite{RST}) give that $v_j+O_j$ and so $O_j$ are planar open disks. Finally, since the Gaussian image of $C_{n,j}$ equals $N_{S_n}(O_{n,j})=B(p_j;\frac1{n})$, then $v_j+O_j$ is contained in a plane orthogonal to $p_j$. 

This proves Claim \ref{cla:discos} and also that the convergence of $\{O_{n,j}\}_{n\geq n_0}$ to $O_j$ is smooth for all $j\in\{1,\ldots,\mgot\}$.

Recall that $\Omega_n=X_n(\Delta_n)$ and set
\[
\widetilde\Omega_n:=\{X_n(p)+N_{S_n}(p)\colon p\in\Omega_n\}\subset \widetilde S_n
\]
the outer parallel surface to $\Omega_n$ at distance $1$. It follows that $\{\widetilde\Omega_n\}_{n\geq n_0}$ converges to the closure of the open subset of $\widetilde S$ consisting of the points at outer distance $1$ from $\Omega\setminus\partial\Omega\subset S$. 
\begin{definition}\label{def:Sigma}
Call $\Sigma$ to the closure of the set of points of $\widetilde S$ at outer distance $1$ from $\Omega\setminus\partial\Omega\subset S$.
\end{definition}

Let us prove that
\begin{claim}\label{cla:Kpart}
$\Sigma$ is a compact $H$-surface (with $H=1/2$) which is $\mathcal C^{1,1}$ up to the boundary.
\end{claim}

Indeed, since $\Omega_n$ is a compact $K$-surface with $K=1$ (see \eqref{eq:f1}), one has that  $\widetilde\Omega_n$ is a compact positively curved $H$-surface with $H=1/2$. In particular, the norm of the second fundamental form of $\widetilde \Omega_n$ is bounded from above by $1$. So, by standard compactness arguments in Constant Mean Curvature Surface Theory (see again Proposition 2.3 and Lemma 2.4 in \cite{RST}), $\Sigma\setminus\partial\Sigma$ is an open $H$-surface and $\{\widetilde\Omega_n\}_{n\geq n_0}$ converges smoothly to $\Sigma\setminus\partial\Sigma$.

Let us show that $\Sigma$ is $\mathcal C^{1,1}$ up to the boundary. Indeed, fix $j\in\{1,\ldots,\mgot\}$ and recall that $O_j$ is a convex domain in $\Pi_j$; see Claim \ref{cla:discos}. Without loss of generality we may assume that $p_j=(0,0,1)$ and $\Pi_j$ is the horizontal plane of equation $x_3=0$. We also assume that $\theta_j< \pi$; the argument for $\theta_j=\pi$ is given in \cite{AS1}.

Let $\widetilde \Pi_j$ be the horizontal plane in $\r^3$ of equation $x_3=r_j^\ast:=\sqrt{1-r_j^2}>0$; obviously $\widetilde\Pi_j=(0,0,r_j^\ast)+\Pi_j$ and $(0,0,r_j^\ast)+O_j\subset\widetilde\Pi_j$.
\begin{definition}\label{def:Dj}
Call $D_j$  the outer parallel convex domain in $\widetilde\Pi_j$ at distance $r_j$ to $(0,0,r_j^\ast )+O_j$. Also set $\gamma_j:=\partial D_j$.
\end{definition}
Observe that $D_j\cap\widetilde S=\emptyset$.

Since $\{\widetilde\Omega_n\}_{n\geq n_0}$ converges smoothly to $\Sigma\setminus\partial\Sigma$ and the Gauss map of $\widetilde \Omega_n$ is a homeomorphism onto $\Delta_n=\s^2\setminus \cup_{j=1}^\mgot B(p_j;r_j+\frac1{n})$ (recall that the Gauss map of $X_n$ is the identity map of $\s^2$), then
\begin{equation}\label{eq:gammaj}
\gamma_j=\Sigma\cap \overline D_j
\end{equation}
and
\begin{equation}\label{eq:K}
M:=\Sigma\cup\bigcup_{j=1}^\mgot D_j\enskip\text{is the boundary surface of a convex body $\Kscr$ in $\r^3$.}
\end{equation}
Notice that $M\cap\widetilde S=\Sigma\setminus\partial\Sigma$. In particular, we emphasize that 
$\Kscr\subsetneq \widetilde\Lscr$.

For each  $k\in \n$ big enough, the intersection  $\alpha_k:= \Sigma\cap\{x_3=r_j^\ast-\frac1{k}\}$  is a smooth $(\Ccal^\infty)$ convex curve in the horizontal plane $\{x_3=r_j^\ast-\frac1{k}\}$ and the sequence $(\alpha_k)_k$ converges to $\gamma_j$ as $k\to \infty.$ We will check that 
the (planar) curvature of $\alpha_k$ is bounded from above independently of $k$, this will imply that  $\gamma_j$ is a $\Ccal^{1,1}$ curve. Since the second fundamental form $II$ of 
$\Sigma$ is bounded this will also imply that $\Sigma$ is $\Ccal^{1,1}$ up to the boundary component $\gamma_j.$ Now, to check the boundedness of the curvature $\kappa$ of $\alpha_k,$ suppose $\alpha_k$ is parameterized by arclength and denote by ${\sf n}$ its  horizontal outward unit normal.  We have the simple relation  $II(\alpha_k^\prime,\alpha_k^\prime)= \kappa \langle {\sf n}, {\sf N}\rangle$. As  $II$ is bounded and $\langle {\sf n}, {\sf N}\rangle$ is uniformly close to $r_j\ne0$, we conclude that $\kappa$ is uniformly bounded from above.

This proves Claim \ref{cla:Kpart}.

Observe also that, up to passing to a subsequence, $\{\Omega_n\}_{n\geq n_0}$ converges smoothly to $\Omega\setminus\partial\Omega$ by the same argument. This also shows that $\Omega\setminus\partial\Omega$ is a $K$-surface (with $K=1$).

We now prove the following
\begin{claim}\label{cla:injective}
$N_{\Sigma\setminus\partial\Sigma}\colon \Sigma\setminus\partial\Sigma\to\s^2$ is a homeomorphism onto $\s^2\setminus\cup_{j=1}^\mgot \overline B(p_j;r_j)$.
\end{claim}
Denote by $N_{\Omega\setminus\partial\Omega}$ and $N_{O_j}$ the outer Gauss map of $\Omega\setminus\partial\Omega$ and $O_j$, $j=1,\ldots,\mgot$, respectively. Since $\Sigma\setminus\partial\Sigma$ is the parallel surface to $\Omega\setminus\partial\Omega$, then to prove Claim \ref{cla:injective} it suffices to show that $N_{\Omega\setminus\partial\Omega}\colon \Omega\setminus\partial\Omega\to\s^2$ is a homeomorphism onto $\s^2\setminus\cup_{j=1}^\mgot \overline B(p_j;r_j)$.

Indeed, let us first check that $N_{\Omega\setminus\partial\Omega}(\Omega\setminus\partial\Omega)\subset \s^2\setminus\cup_{j=1}^\mgot \overline B(p_j;r_j)$ and is injective. Indeed, since $\Lscr$ is convex and $\Omega\setminus\partial\Omega$ is locally strictly convex, then $(T_x \Omega)\cap \Lscr=\{x\}$ for all $x\in \Omega\setminus\partial\Omega$, hence $N_{\Omega\setminus\partial\Omega}$ is injective. (Here $T_x\Omega$ denotes the affine tangent plane to $\Omega$ at $x\in\Omega\setminus\partial\Omega$.) Moreover, since $N_{S_n}(\Omega_n) \subset \s^2\setminus\cup_{j=1}^\mgot \overline B(p_j;r_j)$ for all $n\geq n_0$ and $\Omega_n$ converges smoothly to $\Omega\setminus\partial\Omega$, then $N_{\Omega\setminus\partial\Omega}(\Omega\setminus\partial\Omega)\subset 
\s^2\setminus\cup_{j=1}^\mgot B(p_j;r_j)$. Finally, since $\Omega\setminus\partial\Omega$ is a $K$-surface, then $N_{\Omega\setminus\partial\Omega}$ is an open map and so $N_{\Omega\setminus\partial\Omega}(\Omega\setminus\partial\Omega)\subset \s^2\setminus\cup_{j=1}^\mgot \overline B(p_j;r_j)$. 

Let us now check that $N_{\Omega\setminus\partial\Omega}\colon \Omega\setminus\partial\Omega\to\s^2\setminus\cup_{j=1}^\mgot \overline B(p_j;r_j)$ is surjective. Indeed, take $p\in \s^2\setminus\cup_{j=1}^\mgot \overline B(p_j;r_j)$. Since $\{\Delta_n\}_{n\geq n_0}$ is an exhaustion of $\s^2\setminus\cup_{j=1}^\mgot \overline B_j(p;r_j)$, then there exist $\epsilon>0$ and $n_1\in\n$ such that $B(p;\epsilon)\subset \Delta_n$ for all $n\geq n_1$. Recall that the Gauss map of $X_n$ is the identity map and $\Omega_n=X_n(\Delta_n)$; hence $N_{S_n}(X_n(p))=p$ for all $n\geq n_1$. On the other hand, up to passing to a subsequence, $\{X_n(p)\in \Omega_n\}_{n\geq n_1}$ converges to a point $x\in \Omega\setminus\partial\Omega$ and, since the convergence of $\{\Omega_n\}_{n\geq n_0}$ to $\Omega\setminus\partial\Omega$ is smooth, then $N_{\Omega\setminus\partial\Omega}(x)=\lim_{n\to\infty} N_{S_n}(X_n(p))=p$.

This proves Claim \ref{cla:injective}.

\begin{definition}\label{def:B}
Let $\Bcal$ be the polyhedral region determined by the affine planes $\widetilde\Pi_j$, $j=1,\ldots,\mgot$.
\end{definition}

To finish the proof let us show that the polyhedral region $\Bcal$ and the surface $\Sigma$ satisfy the conclusion of Theorem \ref{th:intro} for $H=1/2$.

Indeed, Claim \ref{cla:discos} implies that $\widetilde\Pi_j$ is orthogonal to $p_j$ for all $j\in\{1,\ldots,\mgot\}$; hence the frontier $\partial\Bcal$ of $\Bcal$ consists of $\mgot$ planar regions $F_1,\ldots,F_\mgot$, where $F_j$ is orthogonal to $p_j$ for all $j\in\{1,\ldots,\mgot\}$. 

By Claim \ref{cla:Kpart}, \eqref{eq:K}, and \eqref{eq:gammaj}, $\Sigma$ is an embedded constant mean curvature surface with $H=1/2$, $\Sigma$ is positively curved, $\Sigma\subset\Bcal$, and $\Sigma\cap\partial\Bcal=\cup_{j=1}^\mgot \gamma_j$. Furthermore $\Sigma\cap F_j=\gamma_j$ and $D_j$ is a convex disk; hence $\gamma_j$ is a convex Jordan curve in $F_j$ for all $j=1,\ldots,\mgot$. This proves conditions {\rm (ii)} and {\rm (iv)} in Theorem \ref{th:intro}.

On the other hand, properties {\rm (i)} and {\rm (iii)} are ensured by Claims \ref{cla:Kpart} and \ref{cla:injective}; take into account \eqref{eq:Bj}. In particular, $\Sigma$ is a capillary surface in $\Bcal$.

To check Theorem \ref{th:intro} {\rm (v)} we argue as follows.  First of all notice that taking limits in \eqref{eq:area} as $n$ goes to infinity, one gets that
\[
{\rm Area}(O_j)=a_j\quad \forall j\in\{1,\ldots,\mgot\}.
\]
(Here we have also taken into account that the convergence of $\{O_{n,j}\}_{n\geq n_0}$ to $O_j$ is smooth; see the proof of Claim \ref{cla:discos}.) Therefore, in view of \eqref{eq:rj}, it suffices to show that
$
{\rm Area}(O_j)={\rm Area}(D_j)-{\rm Length}(\partial D_j) r_j+\pi r_j^2,
$
or equivalently,
\begin{equation}\label{eq:parallel}
{\rm Area}(D_j)-{\rm Area}(D_j')={\rm Length}(\partial D_j) r_j-\pi r_j^2,
\end{equation}
where $D_j':=(0,0,r_j^\ast)+O_j\subset D_j\subset\widetilde\Pi_j$, $j=1,\ldots,\mgot$. {To see this, suppose first that $D_j$ has smooth boundary and consider }the vector field ${\sf I}\colon\widetilde\Pi_j\to\widetilde\Pi_j$ given by ${\sf I}(p)=p$ for all $p\in\widetilde\Pi_j$. Since the divergence of ${\sf I}$ equals $2$ everywhere on $\widetilde\Pi_j$, then the Divergence Theorem gives that
\begin{equation}\label{eq:A2}
2{\rm Area}(D_j\setminus D_j')=\int_{\partial D_j} \langle {\sf I},\eta_{\partial D_j} \rangle + \int_{\partial D_j'} \langle {\sf I},\eta_{\partial D_j'} \rangle,
\end{equation}
where $\eta_{\partial D_j}$ and $\eta_{\partial D_j'}$ are the unit normals to $\partial D_j$ and $\partial D_j'$, respectively, pointing out of $D_j\setminus \overline D_j'$. Since $D_j\subset \widetilde\Pi_j$ is the outer parallel domain at distance $r_j$ to $D_j'$, then 
\begin{equation}\label{eq:A3}
\eta_{\partial D_j'}\big(p-r_j\eta_{\partial D_j}(p)\big)=-\eta_{\partial D_j}(p) \quad \forall p\in \partial D_j.
\end{equation}
Taking into account the Gauss-Bonnet theorem for closed planar curves, \eqref{eq:A3} and a straightforward computation give that
\[
\int_{\partial D_j'} \langle {\sf I},\eta_{\partial D_j'} \rangle=-\int_{\partial D_j} \langle {\sf I},\eta_{\partial D_j} \rangle +2 r_j {\rm Length}(\partial D_j) -2\pi r_j^2.
\]
Combining this equation with \eqref{eq:A2} and the fact that $D_j'\subset D_j$, { one obtains \eqref{eq:parallel} when $\partial D_j$ is smooth. The general case is obtained by applying the previous argument to the smooth curve $(\widetilde\Pi_j-\epsilon p_j)\cap \Sigma$ for $\epsilon >0$ small enough and letting $\epsilon \to 0$, proving condition {\rm (v)}.

This concludes the proof of the existence part in Theorem \ref{th:intro}.}


\section{Uniqueness and the balancing condition}\label{sec:balancing}

In this section we prove that the balancing condition \eqref{eq:balancing} is in fact necessary for the existence of a polyhedral region $\mathcal B$ and a capillary surface $\Sigma$ in $\Bcal$ as in Theorem \ref{th:intro}. We also show that the capillary surfaces furnished by Theorem \ref{th:intro} are unique up to translations.


Let $\{p_1,\ldots,p_\mgot\}\subset\s^2$, let $H>0$, let $(\theta_1,\ldots,\theta_\mgot)\in(\frac{\pi}{2}, \pi)^\mgot$, and let $B_j\subset\s^2$, $j=1,\ldots,\mgot$, be as in the statement of Theorem \ref{th:intro}.

Assume there exists a polyhedral  region $\Bcal$ in $\r^3$ with frontier $\partial\Bcal$ consisting of $\mgot$ planar regions $F_1,\ldots, F_\mgot$, where $F_j$ is orthogonal to $p_j$ for all $j=1,\ldots,\mgot$, and there exists a capillary surface $\Sigma$ which is $\mathcal C^1$ up to the boundary, with constant mean curvature $H$, in $\Bcal$ enjoying properties Theorem \ref{th:intro} {\rm (ii)}--{\rm (iv)}.
Set 
\begin{equation}\label{eq:(v)}
a_j:={\rm Area}(D_j)-{\rm Length}(\partial D_j)\frac{\sin\theta_j}{2H}+\pi\left(\frac{\sin\theta_j}{2H}\right)^2,\quad j=1,\ldots,\mgot.
\end{equation}

Let us prove that
\begin{claim}\label{cla:balancing}
The balancing condition \eqref{eq:balancing} is satisfied. 
\end{claim}
Indeed, up to rescaling, we may assume without loss of generality that $H=1/2$. Since the Gauss map of $\Sigma$ is a diffeomorphism, then the principal curvatures of $\Sigma$ are different from zero. It follows that the parallel surface to $\Sigma\setminus\partial\Sigma$ at signed distance $-1$ is a regular $K$-surface with $K=1$; call $\Omega$ to the closure of that surface. { Observe that $\Omega$ is compact and its  frontier $\partial \Omega$ consists of $\mgot$ components $\alpha_1,\ldots,\alpha_{\mgot}.$ We will see that  each $\alpha_j$ is a convex Jordan curve bounding a convex disk $O_j$ contained in a plane $\Pi_j$ orthogonal to $p_j$. Indeed, fix $j\in\{1,\ldots,\mgot\}$, we may suppose without loss of generality that $p_j=(0,0,1)$ and that the face $F_j$ is contained in the plane 
$\{x_3=0\}$.   
For any $\epsilon>0$ small enough the intersection $\gamma_j^{\epsilon}=\Sigma\cap
\{x_3=-\epsilon\}$ is a smooth ($\Ccal^\infty$) convex curve; denote by $\sf{n}_{\epsilon}$ its outward unit normal in the plane $\{x_3=-\epsilon\}$. Observe that, by {\rm (iii)} and {\rm (iv)}, the Gauss map $N$ of  $\Sigma\setminus\partial\Sigma$ is a diffeomorphism onto $\s^2\setminus\cup_{j=1}^{\mgot} \overline B_j$ and so    $\langle N, {\sf n}_{\epsilon}\rangle 
> \sin\theta_j$ along $\gamma_j^{\epsilon}.$ 
 Parameterize  $\gamma_j^{\epsilon}$ by arclength and denote by $\kappa_{\epsilon}$ its curvature. Let also $II$ be the second fundamental form of $\Sigma$. We have the relation $II({\gamma_j^{\epsilon}}^\prime, {\gamma_j^{\epsilon}}^\prime)
=\kappa_{\epsilon} \langle N, \sf{n}_{\epsilon}\rangle$. As $\Sigma$ is positively curved with mean curvature  $H=1/2$, we have $II({\gamma_j^{\epsilon}}^\prime, {\gamma_j^{\epsilon}}^\prime)<1$ and therefore 
$\kappa_{\epsilon}< 1/\sin\theta_j $. It is a well known fact that this inequality implies that the parallel inward curve to $\gamma_j^{\epsilon}$ at distance $\sin\theta_j$ in the plane $\{x_3=-\epsilon\}$
is a smooth convex curve. Letting $\epsilon \to 0$, we conclude that the parallel inward set, 
$\gamma_j^\ast$, to $\gamma_j$ at distance $\sin\theta_j$ is either a convex curve, a segment, or  a point.  Suppose
$\gamma_j^\ast$ is a segment then $\gamma_j$  is  the boundary of the convex hull of the circles of radius $\sin\theta_j$ centered at its endpoints. In particular  $\Sigma$ is smooth up to $\gamma_j$ (see \cite{dierkes}) and $\gamma_j$  contains two  line segments along each of which the Gauss map is constant. Let $\Gamma$ be one of those segments and let $\mathcal C$ be the cylinder 
of mean curvature $H=1/2$ with a generating line containing $\Gamma$ and with the same mean curvature vector field as $ \Sigma$ along $\Gamma.$  The cylinder $\mathcal C$ clearly exists.
By the maximum principle, $\Sigma$ and $\mathcal C$ coincide in a neighborhood of $\Gamma$ and by analytic extension $\Sigma$ is a part of $\mathcal C.$ This is clearly impossible by the structure of $\gamma_j.$ So $\gamma_j^{\ast}$ cannot be a segment.  Suppose now that  $\gamma_j^{\ast}$ is a point. This means that $\gamma_j$ is a circle of radius $\sin\theta_j.$ Again in this case $\Sigma$ is smooth up to $\gamma_j$ and its Gauss map along $\gamma_j$ makes an angle $\theta_j$ with the horizontal (recall we are assuming $p_j=(0,0,1)$). One can check that there is a (unique) unduloid 
$\mathcal N$ containing $\gamma_j$ and whose mean curvature vector field coincides with that of $\Sigma$ along 
$\gamma_j$. By the maximum principle again we conclude that $\Sigma$ is a part of $\mathcal N$. However  this is excluded by our hypothesis. It is noteworthy that more generally, as  is shown in \cite{GHM}, an embedded isolated singularity of a $K$-surface in $\r^3$ is determined by  the curve 
in $\s^2$ of its limit unit normals at the singularity.

So the sets $\alpha_j$ are convex curves and bound convex domains $O_j\subset \Pi_j$}.
We emphasize that the curves $\alpha_j$ need not be $\Ccal^1$ and might have singularities. Moreover, it follows from \eqref{eq:(v)} and Theorem \ref{th:intro} {\rm (iii)} that ${\rm Area}(O_j)=a_j$ for all $j\in\{1,\ldots,\mgot\}$; use the same argument that proves \eqref{eq:parallel}. 

Notice that $S:=\Omega\cup(\cup_{j=1}^\mgot \overline O_j)$ is the boundary surface of a convex body and the sets $\Omega,O_1,\ldots,O_\mgot$ are pairwise disjoint. Denote by $G_S$ the {\em generalized Gauss map} of $S$, that is, the set-valued map $G_S\colon S\to\s^2$ mapping every point $p\in S$ to the set of all outer normals of the supporting planes of $S$ passing through $p$. We then consider the measure $\mu_S$ induced by $S$ on $\s^2$, which is called the {\em area function} of $S$ and is given by
\[
\mu_S(E)={\rm Area}(\{p\in S\colon G(p)\cap E\neq \emptyset\})\quad \text{for any Borel subset $E\subset\s^2$.}
\]
It trivially follows from the above description of the surface $S$ that
\begin{equation}\label{eq:measure}
\mu_S=\mu_{\s^2}|_{\s^2\setminus\cup_{j=1}^\mgot \overline B_j}+\sum_{j=1}^\mgot a_j \delta_{p_j},
\end{equation}
where $\mu_{\s^2}$ denotes the canonical Lebesgue measure on $\s^2$ and $\delta_{p_j}$ the Dirac measure at the point $p_j$, $j=1,\ldots,\mgot$. As observed by Minkowski, since $S$ is the boundary surface of a convex body then
\begin{equation}\label{eq:genmin1}
\int_{\s^2}{\rm i}_{\s^2}\mu_S=\vec 0,
\end{equation}
where ${\rm i}_{\s^2}\colon\s^2\to\r^3$ is the inclusion map. Denote by ${\bf N}\colon \Omega\setminus\partial\Omega\to\s^2$ and $dA$ the Gauss map and area element of $\Omega\setminus\partial\Omega$, respectively. From \eqref{eq:measure}, 
\begin{equation}\label{eq:genmin2}
\int_{\s^2}{\rm i}_{\s^2}\mu_S=\int_{\Omega\setminus\partial\Omega} {\bf N}\, dA+\sum_{j=1}^\mgot a_j p_j.
\end{equation}
Since ${\bf N}$ is a diffeomorphism onto $\s^2\setminus\cup_{j=1}^\mgot \overline B_j$ and $\Omega\setminus\partial\Omega$ is a $K$-surface with $K=1$, then 
\begin{equation}\label{eq:genmin3}
\int_{\Omega\setminus\partial\Omega} {\bf N}\, dA= \int_{\s^2\setminus\cup_{j=1}^\mgot \overline B_j} p\, dp = 
-\int_{\cup_{j=1}^\mgot B_j} p\, dp = -\sum_{j=1}^\mgot \pi \sin^2\theta_j \, p_j.
\end{equation}
For the latter equality we have taken into account that $B_j=B(p_j;\sin\theta_j)$ and \eqref{eq:spherical1}.

Combining \eqref{eq:genmin1}, \eqref{eq:genmin2}, and \eqref{eq:genmin3} one gets \eqref{eq:balancing}; recall that we are assuming $H=1/2$. 

This proves Claim \ref{cla:balancing}.

Finally, observe that $S$ is the unique up to translations solution to the generalized Minkowski problem for the measure $\mu_S$ \eqref{eq:measure}. See \cite{CY,Sc} (or \cite{AS1}) for a detailed exposition of this classical subject. In particular, we obtain the following result.

\begin{claim}\label{cla:unique}
In the assumptions of Theorem \ref{th:intro}. The polyhedral region $\Bcal$ and the capillary surface $\Sigma$ in $\Bcal$ satisfying properties Theorem \ref{th:intro} {\rm (i)}--{\rm (v)} are unique up to translations.
\end{claim}  
Indeed, just recall that $\Sigma\setminus\partial\Sigma$ is the outer parallel surface to $\Omega\setminus\partial\Omega$ at distance $1$.

This concludes the proof of Theorem \ref{th:intro}.


\section{Proof of Corollary \ref{co:pi2}}\label{sec:pi2}

For the existence part let $\mgot\in\n$, $H>0$, $\theta_\mgot=\frac{\pi}2$, let $\theta_1,\ldots,\theta_{\mgot-1}\in(\frac{\pi}2,\pi)^{\mgot-1}$, and let $p_j$, $B_j$, and $a_j$, $j=1,\ldots,\mgot$, be as in the statement of Theorem \ref{th:intro}.  

Since $\theta_\mgot=\frac{\pi}2$, equation \eqref{eq:balancing} reads
\begin{equation}\label{eq:balancingpi2}
\sum_{j=1}^{\mgot-1} \Big(a_j-\frac{\pi}{4H^2} \sin^2\theta_j \Big) p_j =\Big(\frac{\pi}{4H^2} -a_\mgot \Big) p_\mgot,
\end{equation}
and $B_\mgot=\{q\in\s^2\colon \langle q, p_\mgot\rangle>0\}$. Since the closed disks $\overline B_1,\ldots,\overline B_\mgot$ are pairwise disjoint, then 
\begin{equation}\label{eq:Bpi2}
\overline B_j\subset \s^2\setminus \overline B_\mgot=\{q\in\s^2\colon \langle q, p_\mgot\rangle<0\}\quad \text{for all $j=1,\ldots,\mgot-1$.}
\end{equation}

Denote by $s\colon\r^3\to\r^3$ the symmetry with respect to the vectorial plane $\Pi$ in $\r^3$ orthogonal to $p_\mgot$. Observe that
\begin{equation}\label{eq:spm}
s(p_\mgot)=-p_\mgot,
\end{equation}
\begin{equation}\label{eq:sBm}
s(\s^2\setminus \overline B_\mgot)=B_\mgot, \quad\text{and}\quad \Pi\cap \overline B_j=\emptyset\;\; \forall j=1,\ldots,\mgot-1.
\end{equation}
Label $p_j'=p_j$, $\theta_j'=\theta_j$, $B_j'=B_j$, and $a_j'= a_j$, and set $p_{\mgot-1+j}'= s(p_j')$, $\theta_{\mgot-1+j}'=\theta_j'$, $B_{\mgot-1+j}'= s(B_j)$, and $a_{\mgot-1+j}'=a_j'$, for all $j=1,\ldots,\mgot-1$. In view of \eqref{eq:Bpi2} and \eqref{eq:sBm}, the closed disks $\overline B_j'$, $j=1,\ldots,2\mgot-2$, are pairwise disjoint. Moreover,
\begin{eqnarray*}
\sum_{j=1}^{2\mgot-2} \Big(a_j'-\frac{\pi}{4H^2} \sin^2\theta_j' \Big) p_j' & = & \sum_{j=1}^{\mgot-1} \Big(a_j-\frac{\pi}{4H^2} \sin^2\theta_j \Big) (p_j+ s(p_j))\\
 & \stackrel{\text{\eqref{eq:balancingpi2}}}{=} & \Big(\frac{\pi}{4H^2}-a_\mgot\Big) (p_\mgot+ s(p_\mgot)) \stackrel{\text{\eqref{eq:spm}}}{=} 0.
\end{eqnarray*}
Therefore Theorem \ref{th:intro} applies to the data $p_j'$, $\theta_j'$, $a_j'$, $j=1,\ldots,2\mgot-2$ (recall that $\theta_j'\in(\frac{\pi}2,\pi)$ for all $j$). This ensures the existence of a unique up to translations polyhedral region $\Bcal'$ in $\r^3$ with frontier $\partial \Bcal'$ consisting of $2\mgot-2$ planar regions $F_1',\ldots,F_{2\mgot-2}'$, where $F_j'$ is orthogonal to $p_j'$ for all $j\in\{1,\ldots,2\mgot-2\}$, and a unique up to translations embedded capillary surface $\Sigma'$ in $\Bcal'$, with constant mean curvature $H$, satisfying the following properties:
\begin{enumerate}[\rm (i$'$)]
\item $\Sigma'$ is $\Ccal^{1,1}$ up to the boundary.
\item $\gamma_j':=\Sigma'\cap F_j'$ is a convex Jordan curve contained in the relative interior of $F_j'$. 
\item $\Sigma'$ meets $F_j'$ at constant angle $\theta_j'$ along $\gamma_j'$ for all $j\in\{1,\ldots,2\mgot-2\}$.
\item $\Sigma'\setminus\partial\Sigma'$ is positively curved and $\Sigma'\cup (\cup_{j=1}^\mgot D_j')$ is the boundary surface of a convex body in $\r^3,$ where $D_j'$ is the convex connected component of  $F_j'\setminus\gamma_j'$ for all $j\in\{1,\ldots,2\mgot-2\}$.
\item $a_j'= {\rm Area}(D_j')-{\rm Length}(\partial D_j')\frac{\sin\theta_j'}{2H} + \pi \left(\frac{\sin\theta_j'}{2H}\right)^2$  for all $ j\in\{1,\ldots,2\mgot-2\}.$
\end{enumerate}

From the definition of the data $p_j'$, $\theta_j'$, and $a_j'$, $j=1,\ldots,2\mgot-2$, we infer that the polyhedral region $s(\Bcal')$ and the capillary surface $s(\Sigma')$ in $s(\Bcal')$ also meet the above conditions; hence, up to a translation, $s(\Bcal')=\Bcal'$ and $s(\Sigma')=\Sigma'$. This simply means that $\Bcal'$ and $\Sigma'$ are symmetric with respect to $\Pi$. Set $\Bcal:=\Bcal'\cap\Pi_0$ and $\Sigma:=\Sigma'\cap\Pi_0$ where $\Pi_0$ is the closure of the connected component of $\r^3\setminus\Pi$ containing $-p_\mgot$.

Let us check that the polyhedral region $\Bcal$ and the capillary surface $\Sigma$ in $\Bcal$ satisfy the conclusion of Theorem \ref{th:intro}. Indeed, observe that the faces of $\Bcal$ are $F_j:=F_j'$, $j=1,\ldots,\mgot-1$, and $F_\mgot:=\Bcal\cap\Pi$; hence $F_j$ is orthogonal to $p_j$ for all $j=1,\ldots,\mgot$. Since $\Sigma=\Sigma'\cap\Pi_0\subset\Sigma'$, then {\rm (i$'$)} and {\rm (iv$'$)}, respectively, imply properties {\rm (i)} and {\rm (iv)} in Theorem \ref{th:intro}. In view of {\rm (ii$'$)}, {\rm (iii$'$)}, and {\rm (v$'$)}, it suffices to check properties {\rm (ii)}, {\rm (iii)}, and {\rm (v)} for $j=\mgot$. Indeed, since $\gamma_\mgot=\Sigma'\cap\Pi\subset\Bcal'\setminus\partial \Bcal'$, then {\rm (iv$'$)} implies {\rm (ii)} for $j=\mgot$. By the symmetry, $\Sigma$ intersects $\Pi$ at constant angle $\frac{\pi}2$, proving {\rm (iii)}. On the other hand, the same argument that proved equation \eqref{eq:parallel} in Sec.\ \ref{sec:proof} now gives {\rm (v)} for $j=\mgot$. Finally, let $\Bcal_0$ be a polyhedral region and let $\Sigma_0$ be a capillary surface in $\Bcal_0$ satisfying properties {\rm (i)}--{\rm (v)}. Translate $\Bcal_0$ and $\Sigma_0$ so that $\Bcal_0\cap\Pi$ is the face of $\Bcal_0$ orthogonal to $p_\mgot$. Since $\theta_\mgot=\frac{\pi}2$, then $\Bcal_0\cup s(\Bcal_0)$ is a polyhedral region and $\Sigma_0 \cup s(\Sigma_0)$ is a capillary surface in $\Bcal_0 \cup s(\Bcal_0)$ which satisfy conditions {\rm (i$'$)}--{\rm (v$'$)}. Therefore the uniqueness of Theorem \ref{th:intro} ensures that, up to a translation, $\Bcal_0\cup s(\Bcal_0)=\Bcal'$ and $\Sigma_0 \cup s(\Sigma_0)=\Sigma'$, hence $\Bcal_0=\Bcal$ and $\Sigma_0=\Sigma$, i.e. $\Bcal$ and $\Sigma$ are unique up to translations. 

 For the converse part of Theorem \ref{th:intro} consider  a data $(\mathcal B, \Sigma)$ satisfying the required hypotheses with $p_\mgot=\frac{\pi}{2}$. Call $\Pi$ the affine plane containing the face
 $F_\mgot $ and $s\colon\r^3\to\r^3$ the symmetry with respect to  $\Pi$.  Since $\Sigma$ is orthogonal to $\Pi$ and has constant mean curvature, the surface 
  $\Sigma' =\Sigma\cup s(\Sigma)$ is regular along $\Sigma\cap\Pi$ and we can apply the argument  in Sec.\ \ref{sec:balancing} to the data $(\mathcal B^\prime,\Sigma^\prime)$ where 
  $\mathcal B'=\mathcal B \cup s(\mathcal B)$. It follows that the parallel surface to $\Sigma'$ at distance $-1$ (assuming again after rescaling that $H=1/2$) is an embedded $K$-surface 
  bounded by convex curves contained in affine planes. In particular the parallel surface at distance $-1$ to $\Sigma$ has the same property and the rest of the proof of Claim \ref{cla:balancing} applies to $\Sigma$. 


\section{Final comments}\label{sec:comments}

It would be interesting to decide whether the capillary surfaces we have obtained enjoy higher boundary regularity. In this respect, we point out that Choe \cite{choe} has proven that a minimal surface in $\r^3$ meeting a plane along a $\mathcal C^1$ arc with a constant angle can be extended analytically across the arc, and raised the question of extending his result to the case of constant mean curvature surfaces. 

Finally, we also think it is an  interesting problem to decide whether these surfaces are minimizing for the energy functional (\ref{eq:energy}) or at least whether they are stable. We refer to \cite{ros} and \cite{ros-souam} for the discussion of the  notion of stability in the capillarity setting  and for some results  on this issue
in the free boundary case inside mean-convex domains and in the general capillary case inside a ball, respectively.


\end{document}